\documentclass[a4paper,12pt]{amsart}
\usepackage{amssymb}
\usepackage{ifthen}
\usepackage[dvips]{graphicx}
\nonstopmode \numberwithin{equation}{section}
\usepackage{amssymb}
\usepackage{ifthen}
\usepackage{graphicx}
\usepackage{amsmath}
\usepackage[T1]{fontenc} %skandit
\usepackage[a4paper,textwidth=15cm,textheight=25cm,left=4cm]{geometry}
\nonstopmode \numberwithin{equation}{section}
\theoremstyle{plain}
\newtheorem{thm}[equation]{Theorem}
\newtheorem{cor}[equation]{Corollary}
\newtheorem{lem}[equation]{Lemma}
\newtheorem{prop}{Proposition}

\newtheorem{conj}{Conjecture}
\newtheorem{innercustomthm}{Theorem}
\newenvironment{customthm}[1]
{\renewcommand\theinnercustomthm{#1}\innercustomthm}
{\endinnercustomthm}

\theoremstyle{definition}
\newtheorem{defn}{Definition}[section]

\newtheorem{prob}{Problem}
\newtheorem{rem}{Remark}[section]

\newcounter{minutes}
\divide\time by 60
\newcounter{hours}
\multiply\time by 60
\addtocounter{minutes}{-\time}
\newcounter {own}
\def\theown {\thesection       .\arabic{own}}

\newenvironment{pf}[1][]{%
	\vskip 3mm
	\noindent
	\ifthenelse{\equal{#1}{}}%
	{{\slshape Proof. }}%
	{{\slshape #1.} }%
}%
{\qed\bigskip}

\newcounter{alphabet}
\newcounter{tmp}

\newcommand{\real}{{\operatorname{Re}\,}}

\def\be{\begin{equation}}
\def\ee{\end{equation}}

\newcommand{\bee}{\begin{enumerate}}
	\newcommand{\eee}{\end{enumerate}}

\newcommand{\blem}{\begin{lem}}
	\newcommand{\elem}{\end{lem}}
\newcommand{\bthm}{\begin{thm}}
	\newcommand{\ethm}{\end{thm}}
\newcommand{\bcor}{\begin{cor}}
	\newcommand{\ecor}{\end{cor}}
\newcommand{\beg}{\begin{examp}}
	\newcommand{\eeg}{\end{examp}}
\newcommand{\begs}{\begin{examples}}
	\newcommand{\eegs}{\end{examples}}

\newcommand{\bdefn}{\begin{defn}}
	\newcommand{\edefn}{\end{defn}}

\newcommand{\bprob}{\begin{prob}}
	\newcommand{\eprob}{\end{prob}}
\newcommand{\bei}{\begin{itemize}}
	\newcommand{\eei}{\end{itemize}}

\newcommand{\bcon}{\begin{conj}}
	\newcommand{\econ}{\end{conj}}
\newcommand{\bcons}{\begin{conjs}}
	\newcommand{\econs}{\end{conjs}}
\newcommand{\bprop}{\begin{prop}}
	\newcommand{\eprop}{\end{prop}}
\newcommand{\br}{\begin{rem}}
	\newcommand{\er}{\end{rem}}
\newcommand{\brs}{\begin{rems}}
	\newcommand{\ers}{\end{rems}}
\newcommand{\bo}{\begin{obser}}
	\newcommand{\eo}{\end{obser}}
\newcommand{\bos}{\begin{obsers}}
	\newcommand{\eos}{\end{obsers}}
\newcommand{\bpf}{\begin{pf}}
	\newcommand{\epf}{\end{pf}}
\newcommand{\ba}{\begin{array}}
	\newcommand{\ea}{\end{array}}
\newcommand{\beq}{\begin{eqnarray}}
\newcommand{\beqq}{\begin{eqnarray*}}
\newcommand{\eeq}{\end{eqnarray}}
\newcommand{\eeqq}{\end{eqnarray*}}

\begin{document}

\title{Bohr phenomenon for certain close-to-convex analytic functions}

\author{Vasudevarao Allu}
\address{Vasudevarao Allu,
	School of Basic Science,
	Indian Institute of Technology Bhubaneswar,
	Bhubaneswar-752050, Odisha, India.}
\email{avrao@iitbbs.ac.in}

\author{Himadri Halder}
\address{Himadri Halder,
	School of Basic Science,
	Indian Institute of Technology Bhubaneswar,
	Bhubaneswar-752050, Odisha, India.}
\email{hh11@iitbbs.ac.in}

\subjclass[{AMS} Subject Classification:]{Primary 30C45, 30C50, 30C80}

\def\thefootnote{}
\footnotetext{ {\tiny File:~\jobname.tex,
		printed: \number\year-\number\month-\number\day,
		\thehours.\ifnum\theminutes<10{0}\fi\theminutes }
} \makeatletter\def\thefootnote{\@arabic\c@footnote}\makeatother

\begin{abstract}
We say that a class $\mathcal{B}$ of analytic functions $f$ of the form $f(z)=\sum_{n=0}^{\infty} a_{n}z^{n}$ in the unit disk $\mathbb{D}:=\{z\in \mathbb{C}: |z|<1\}$ satisfies a Bohr phenomenon if for the largest radius $R_{f}<1$, the following inequality 
$$
\sum\limits_{n=1}^{\infty} |a_{n}z^{n}| \leq d(f(0),\partial f(\mathbb{D}) ) 
$$
holds for $|z|=r\leq R_{f}$ and for all functions $f \in \mathcal{B}$. The largest radius $R_{f}$ is called Bohr radius for the class $\mathcal{B}$.
In this article, we obtain Bohr radius for certain subclasses of close-to-convex analytic functions. We establish the Bohr phenomenon for certain analytic classes $\mathcal{S}_{c}^{*}(\phi),\,\mathcal{C}_{c}(\phi),\, \mathcal{C}_{s}^{*}(\phi),\, \mathcal{K}_{s}(\phi)$. Using Bohr phenomenon for subordination classes \cite[Lemma 1]{bhowmik-2018}, we obtain some radius $R_{f}$ such that Bohr phenomenon for these classes holds for $|z|=r\leq R_{f}$. Generally, in this case  $R_{f}$ need not be sharp, but we show that under some additional conditions on $\phi$, the radius $R_{f}$ becomes sharp bound. As a consequence of these results, we obtain several interesting corollaries on Bohr phenomenon for the aforesaid classes.
\end{abstract}

\maketitle
\pagestyle{myheadings}
\markboth{Vasudevarao Allu and  Himadri Halder}{Bohr phenomenon for certain close-to-convex analytic functions}

\section{Introduction and Preliminaries}
Let $f$ be an analytic function in the unit disk $\mathbb{D}:=\{z\in \mathbb{C}: |z|<1\}$ with the following power series representation
\begin{equation} \label{him-p3-e-1.1}
f(z)=\sum_{n=0}^{\infty} a_{n}z^{n}.
\end{equation}
Then the majorant series $M_{f}(r)$ associated with $f$ given by \eqref{him-p3-e-1.1}, is defined by  
$
M_{f}(r):= \sum_{n=0}^{\infty} |a_{n}|r^{n}$ for $|z|=r<1$. The classical result of H. Bohr \cite{Bohr-1914}, which in the sharp form has been independently proved by Weiner, Riesz and Schur reads as follows:
\begin{customthm}{A}
	Let $f$ be analytic in $\mathbb{D}$ of the form \eqref{him-p3-e-1.1} and $|f(z)|<1$ for all $z\in \mathbb{D}$. Then the associated majorant series 
	\begin{equation} \label{him-p3-e-1.2}
	M_{f}(r)= \sum_{n=0}^{\infty} |a_{n}|r^{n} \leq 1 \quad \mbox{for} \quad |z|=r \leq 1/3
	\end{equation}
	and the constant $1/3$, referred to as the Bohr radius, cannot be improved.
\end{customthm}
In the recent years, studying the Bohr radus has become an interesting problem in various directions in functions of one and several complex variables. The notion of Bohr radius has been extended to several complex variables, to planar harmonic mappings, to polynomials, to solutions of elliptic partial differential equations, and to more abstract settings. For more information and intriguing aspects about Bohr radius and Bohr inequality as stated above, we suggest the reader to glance through the articles \cite{aizn-2000,alkhaleefah-2019,Himadri-Vasu-P2,boas-1997} and the references therein. 

\vspace{3mm}
The inequality \eqref{him-p3-e-1.2} can also be written in the following form 
\begin{equation} \label{him-p3-e-1.3}
\sum_{n=1}^{\infty} |a_{n}z^{n}|\leq 1-|a_{0}|=d(f(0),\partial f(\mathbb{D}))
\end{equation}
for $|z|=r \leq 1/3$, where $d$ is the Euclidean distance. It is worth noting that the existence of the radius $1/3$ in \eqref{him-p3-e-1.3} is independent of the coefficients of the power series \eqref{him-p3-e-1.1}. Analytic functions of the form \eqref{him-p3-e-1.1} with modulus less than $1$ satisfying the inequality \eqref{him-p3-e-1.3}, are sometimes  said to satisfy the classical Bohr phenomenon. Therefore we conclude that Bohr phenomenon occurs in the class of analytic self-maps of the unit disk $\mathbb{D}$. The notion of Bohr phenomenon has been extended to the class of analytic functions from $\mathbb{D}$ into a given domain $D \subseteq \mathbb{C}$. Let $\mathcal{G}$ be the class of analytic functions of the form \eqref{him-p3-e-1.1} which map $\mathbb{D}$ into a given domain $D$ such that $f(\mathbb{D}) \subseteq D$. Suppose there exists the largest radius $r_{D}>0$ such that 
\begin{equation} \label{him-p3-e-1.4}
\sum\limits_{n=1}^{\infty} |a_{n}z^{n}| \leq d(f(0),\partial f(\mathbb{D}) ) \quad \mbox{in} \quad |z| \leq r_{D}
\end{equation}
for all functions $f \in \mathcal{G}$. In this case, we say that $\mathcal{G}$ satisfies the Bohr phenomenon. It has been proved \cite{aizn-2007} that the largest radius $r_{D}$ for convex domain $D$ coincides with the classical Bohr radius $1/3$ while Abu-Muhanna \cite{Abu-2010} has obtained $r_{D}=3-2\sqrt{2}$ for any proper simply connected domain $D$. For more intriguing aspects of Bohr phenomenon, we refer the reader to the articles  \cite{abu-2011,abu-2013,Ali-2017,Ali-2019}. The Bohr phenomenon for certain subclasses of harmonic mappings has also been extensively studied by several authors \cite{abu-2014,Himadri-Vasu-P1,kayumov-2018-b}. 

\vspace{4mm}
Let $\mathcal{A}$ denote the class of normalized analytic functions in $\mathbb{D}$ of the form 
\begin{equation} \label{him-p3-e-1.5}
f(z)=z+ \sum_{n=2}^{\infty} a_{n}z^{n}
\end{equation}
and $\mathcal{S}$ be its standard subclass made up of normalized univalent ({\it i.e.} one-to-one) functions in $\mathbb{D}$. 
%Let $\mathcal{S}^{*}$ (respectively $\mathcal{C}$) be the subclass of $\mathcal{S}$ consisting of starlike (respectively convex) functions in $\mathbb{D}$. It is well-known that  $f \in \mathcal{S}^{*}$ ($\mathcal{C}$ respectively) if, and only, if 
%$\real({zf'(z)}/{f(z)})>0$  for $z\in\mathbb{D}$   $(\real(1+{zf''(z)}/{f'(z)})>0$ for $z\in\mathbb{D}$ respectively). Let $\mathcal{S}^{*}(\alpha) \left(\mathcal{C}(\alpha) \, \mbox{respectively} \right)$ denote the class of all functions in $\mathcal{S}$ such that $\real({zf'(z)}/{f(z)})>\alpha \, \left(\real(1+{zf''(z)}/{f'(z)})>\alpha \,\,
 %\mbox{respectively} \right)$  for $z\in\mathbb{D}$. Clearly, $f \in \mathcal{C}(\alpha)$ if, and only, if $zf' \in \mathcal{S}^{*}(\alpha)$. Note that for $\alpha=0$ the classes $\mathcal{S}^{*}(0):=\mathcal{S}^{*}$ and $\mathcal{C}(0):=\mathcal{C}$ are the well-known classes of starlike and convex functions in $\mathbb{D}$ respectively.  
An analytic function $f$ in $\mathbb{D}$ is said to be subordinate to an analytic function $g$ in $\mathbb{D}$, denoted by $f \prec g$ (sometimes written as $f(z) \prec g(z)$), if $f(z)=g(\omega(z))$ for $z \in \mathbb{D}$, where $\omega : \mathbb{D} \rightarrow \mathbb{D}$ is an analytic function such that $\omega(0)=0$. In particular, when $g$ is univalent in $\mathbb{D}$, then $f\prec g$ if, and only if, $f(0)=g(0)$ and $f(\mathbb{D}) \subseteq g(\mathbb{D})$. Let $\phi : \mathbb{D} \rightarrow \mathbb{C}$ be Ma-Minda function which is analytic and univalent in $\mathbb{D}$ such that $\phi(\mathbb{D})$ has positive real part, symmetric with respect to the real axis, starlike with respect to $\phi(0)=1$ and $\phi ' (0)>0$. Such Ma-Minda functions have the series representation of the form $\phi(z)=1+ \sum_{n=1}^{\infty} B_{n}z^{n} \quad (B_{1}>0)$. For such $\phi$, Ma-Minda \cite{ma minda-1992-a} have considered the classes $\mathcal{S}^{*}(\phi)$ and $\mathcal{C}(\phi)$, called Ma-Minda type starlike and Ma-Minda type convex classes associated with $\phi$ respectively, where  $\mathcal{S}^{*}(\phi)$ and $\mathcal{C}(\phi)$ are the subclasses of functions in $\mathcal{S}$ such that 
$zf'(z)/f(z) \prec \phi (z) \quad \mbox{and} \quad 1+ zf''(z)/f'(z) \prec \phi (z)$
respectively. Clearly, $f \in \mathcal{C}(\phi)$ if, and only if, $zf' \in \mathcal{S}^{*}(\phi)$. It is important to note that for every such $\phi$ described as above, $\mathcal{S}^{*}(\phi) \left( \mathcal{C}(\phi) \, \, \mbox{respectively}\right)$ always a subclass of the well-known starlike class $\mathcal{S}^{*} \left(\mbox{convex class}\,\, \mathcal{C}\, \, \mbox{respectively} \right)$ by taking $\phi(z)=(1+z)/(1-z)$. For more intriguing aspects and geometric properties of starlike and convex functions, we refer the book \cite{vasu-book}. For various $\phi$, the classes $\mathcal{S}^{*}(\phi)$ and $\mathcal{C}(\phi)$ yield various important subclasses of starlike and convex functions, respectively. When $\phi(z)=(1+(1-2\alpha))/(1-z)$, we obtain the classes $\mathcal{S}^{*}(\alpha)$ and $\mathcal{C}(\alpha)$. By taking  $\phi(z)=(1+Az)/(1+Bz)$, $\mathcal{S}^{*}(\phi)$ and $\mathcal{C}(\phi)$ reduce to the Janowski starlike class $\mathcal{S}^{*}[A,B]$ and Janowski convex class $\mathcal{C}[A,B]$ respectively. By taking $\phi(z)=\left((1+z)/(1-z)\right)^{\alpha}$ for $0< \alpha \leq 1$, we obtain the classes of strongly convex and strongly starlike functions of order $\alpha$. By choosing $\phi(z)=(1+sz)^{2}$ with $0< s \leq 1/\sqrt{2}$, the class $\mathcal{S}^{*}(\phi)$ reduces to $\mathcal{ST}_{L}(s):=\mathcal{S}^{*}\left(\left(1+sz\right)^{2}\right)$. Masih and Kanas \cite{masih-2020} have considered the class $\mathcal{ST}_{L}(s)$. Khatter {\it et al.} \cite{khatter-2019} have introduced the class $\mathcal{S}^{*}_{\alpha,e}:=\mathcal{S}^{*}(\alpha +(1-\alpha)e^{z})$ for $0 \leq \alpha <1$.

\vspace{3mm}
The extremal functions $k$ and $h$ respectively for the classes $\mathcal{C}(\phi)$ and $\mathcal{S}^{*}(\phi)$ as follows:
\begin{equation} \label{him-p3-e-1.6}
1+ \dfrac{zk''(z)}{k'(z)} = \phi (z) \quad \mbox{and} \quad \dfrac{zh'(z)}{h(z)} =\phi (z)
\end{equation}
with the normalizations $k(0)=k'(0)-1=0$ and $h(0)=h'(0)-1=0$. The  functions $k$ and $h$ belong to the classes $\mathcal{C}(\phi)$ and $\mathcal{S}^{*}(\phi)$ and they play the role of Koebe functions in the respective classes.
  Ma and Minda \cite{ma minda-1992-a} have obtained the following subordination result and growth estimates for the classes $\mathcal{S}^{*}(\phi)$ and  $\mathcal{C}(\phi)$.
\begin{lem} \label{him-p3-lem-1.6} \cite{ma minda-1992-a}
	Let $f \in \mathcal{S}^{*}(\phi)$. Then $zf'(z)/f(z) \prec zh'(z)/h(z)$ and $f(z)/z \prec h(z)/z$.
\end{lem}
%\begin{lem} \label{him-p3-lem-1.7} \cite{ma minda-1992-a}
%	Assume $f \in \mathcal{S}^{*}(\phi)$ and $|z|=r<1$. Then 
%	\begin{equation} \label{him-p3-e-1.8}
%	-h(-r) \leq |f(z)| \leq h(r).
%	\end{equation} 
%	Equality holds for some $z \neq 0$ if, and only, if f is a rotation of $h$.
%\end{lem}

\begin{lem} \label{him-p2-lem-1.13} \cite{ma minda-1992-a}
	Let $f \in \mathcal{C}(\phi)$. Then $zf''(z)/f'(z) \prec zk''(z)/k'(z)$ and $f'(z) \prec k'(z)$.
\end{lem}

%\begin{lem} \label{him-p2-lem-1.14} \cite{ma minda-1992-a}
%	Assume $f \in \mathcal{C}(\phi)$ and $|z|=r<1$. Then 
%	\begin{equation} \label{him-p2-e-1.15}
%	-k(-r) \leq |f(z)| \leq k(r).
%	\end{equation}
%	Equality holds for some $z \neq 0$ if, and only, if f is a rotation of $k$.
%\end{lem}

%\begin{lem}
%If $f(z)=z+a_{k+1}z^{k+1}+\cdots \in \mathcal{C}(\phi)$, then 
%\end{lem}
Ma-Minda functions $\phi$ have been considered with the condition $\phi '(0)>0$. Motivated by this, recently, Kumar and Banga \cite{sivaprasadkumar-2020} have introduced the function $\Phi$, called non-Ma-Minda function, with the condition $\Phi'(0)<0$ and the other conditions on $\Phi$ are same as that of $\phi$. Note that $\Phi$ can obtained from $\phi$ by a rotation, namely, $z$ by $-z$. 
 By going a similar manner as the definition of $\mathcal{S}^{*}(\phi)$ and $\mathcal{C}(\phi)$ (see \cite{ma minda-1992-a}), Kumar and Banga have considered the classes $\mathcal{S}^{*}(\Phi)$ and $\mathcal{C}(\Phi)$ and studied the growth estimates and other basic properties of these classes.

\vspace{4mm}
 A function $f \in \mathcal{A}$ is said to be close-to-convex if there exists $g \in \mathcal{S}^{*}$ such that $\real \left(zf'(z)/g(z)\right)>0$ for $z\in \mathbb{D}$. Let $\mathcal{K}$ denote the class of close-to-convex functions in $\mathbb{D}$. In 1959, Sakaguchi \cite{sakaguchi-1959} introduced the subclass $\mathcal{S}^{*}_{s}$ of functions starlike with respect to symmetric points, which consists of functions $ f \in \mathcal{S}$ satisfying the condition $\real \left(zf'(z)/\left(f(z)-f(-z)\right)\right)>0$ for $z \in \mathbb{D}$. Motivated by $\mathcal{S}^{*}_{s}$, Wang {\it et.al.} \cite{Wang-2006-b} have considered the class $\mathcal{C}_{s}$. More precisely, a function  $f \in \mathcal{C}_{s}$ if $f$ satisfies the inequality $ 
\real\left((zf'(z))'/\left(\left(f(z)-f(-z)\right)'\right)\right)>0$ in $\mathbb{D}$. A function $f \in \mathcal{A}$ is starlike with respect to conjugate points and convex with respect to conjugate points in $\mathbb{D}$ respectively if $f$ satisfies the conditions 
$$\real \left( \dfrac{zf'(z)}{f(z)+\overline{f(\bar{z})}}\right)>0 \quad \mbox{and} \quad  \real \left( \dfrac{(zf'(z))'}{\left(f(z)+\overline{f(\bar{z})}\right)'}\right)>0 \quad \mbox{for} \quad z \in \mathbb{D}
$$ 
respectively.
A function $f \in \mathcal{A}$ is starlike with respect to symmetric conjugate points in $\mathbb{D}$ if it  satisfies the inequality 
$$
\real \left( \dfrac{zf'(z)}{f(z)-\overline{f(-\bar{z})}}\right)>0, \, z \in \mathbb{D}.
$$
In more general, Ravichandran \cite{Ravichandran-2004} has defined the classes $\mathcal{S}^{*}_{s}(\phi)$ and $\mathcal{C}_{s}(\phi)$. 
\begin{defn}  \cite{Ravichandran-2004} \label{him-p3-def-1.1}
A function $f \in \mathcal{A}$ is in the class $\mathcal{S}^{*}_{s}(\phi)$ if 
$$
\dfrac{2zf'(z)}{f(z)-f(-z)} \prec \phi (z), \quad z\in \mathbb{D}
$$ 
and is in the class $\mathcal{C}_{s}(\phi)$ if 
$$
\dfrac{2(zf'(z))'}{f'(z)+f'(-z)} \prec \phi (z), \quad \quad z\in \mathbb{D}.
$$
\end{defn} 

Similarly, let $\mathcal{S}^{*}_{c}(\phi)$ and $\mathcal{S}^{*}_{sc}(\phi)$ be the corresponding classes of starlike functions with respect to conjugate points and symmetric conjugate points respectively. Let $\mathcal{C}_{c}(\phi)$ and $\mathcal{C}_{sc}(\phi)$ be the corresponding classes of convex functions with respect to conjugate points and symmetric conjugate points respectively
The following lemmas are required to prove our main results.
\begin{lem} \cite{Ravichandran-2004} \label{him-p3-lem-1.12}
Let $\min _{|z|=r} |\phi (z)|=\phi (-r)$, $\max _{|z|=r}|\phi (z)|=\phi (r)$, $|z|=r$. If $f \in \mathcal{C}_{s}(\phi)$, then 
$$
\frac{1}{r} \int \limits _{0}^{r} \phi (-r) (k'(-r^{2}))^{1/2} \, dr\, \leq \, |f'(z)| \, \leq \frac{1}{r} \int \limits _{0}^{r} \phi (r) (k'(r^{2}))^{1/2} \, dr.
$$
\end{lem}
From \cite[Theorem 9]{Wang-2006-b}, for $f \in \mathcal{C}_{s}(\phi)$, we have 
\begin{equation} \label{him-p3-e-13-a}
\int_{0}^{r} \frac{1}{s} \int \limits _{0}^{s} \phi (-t) (k'(-t^{2}))^{1/2} \, dt\,ds\, \leq |f(z)| \leq 
\int_{0}^{r} \frac{1}{s} \int \limits _{0}^{s} \phi (t) (k'(t^{2}))^{1/2} \, dt \, ds
\end{equation}
and the results are sharp for the following function 
\begin{equation} \label{him-p3-e-13-b}
f(z)=\int_{0}^{z} \frac{1}{\xi} \int \limits _{0}^{\xi} \phi (-\eta) (k'(-\eta^{2}))^{1/2} \, d\eta \, d\xi \, \in \mathcal{C}_{s}(\phi),
\end{equation}
since it belongs to the class $\mathcal{C}(\phi)$ and having real coefficients.

\begin{lem} \label{him-p3-lem-1.16} \cite{graham-1996}
Let $f(z)=z+a_{l+1}z^{l+1}+\cdots \in \mathcal{C}(\phi)$, then we have 
$$
(k'(-r^{l}))^{1/l}\leq |f'(z)|\leq (k'(r^{l}))^{1/l}.
$$
The bounds are sharp for some suitable rotations of the function $K_{l}$ which is defined by
$$
K_{l}(z)=\int \limits  _{0}^{z}(k'(\xi^{l}))^{1/l}\, d\xi, \quad z\in \mathbb{D},
$$
where $k$ is defined in \eqref{him-p3-e-1.6}.
\end{lem}
In particular for $l=2$ we can obtain the bounds of $|f'(z)|$ for odd convex functions. From Lemma \ref{him-p3-lem-1.16}, the following can be easily obtained for $l=2$
$$
\int \limits _{0}^{r} (k'(-t^{2}))^{1/2}\, dt \leq |f(z)| \leq \int \limits_{0}^{r} (k'(t^{2}))^{1/2}\, dt.
$$
The result is sharp for the function $K_{2}=K$ is defined by  $K(z):= \int  _{0}^{z}(k'(\xi^{2}))^{1/2}\, d\xi$. It is easy to see that $K$ is odd convex function which belongs to $\mathcal{C}(\phi)$. Similary, we note that the function $H$ is defined by  $ H(z):=(h(z^{2}))^{1/2}$ is a Koebe type function for odd starlike class in $\mathcal{S}^{*}(\phi)$ and satisfies the relation
 \begin{equation} \label{him-p3-e-1.17}
 zK'(z)=H(z).
 \end{equation}
  
\begin{lem} \cite{Ravichandran-2004} \label{him-p3-lem-1.13}
Let $\min _{|z|=r} |\phi (z)|=\phi (-r)$, $\max _{|z|=r}|\phi (z)|=\phi (r)$, $|z|=r$. If $f \in \mathcal{S}^{*}_{c}(\phi)$, then 
\begin{enumerate}
	\item [(i)] $h'(-r) \leq |f'(z)| \leq h'(r)$ \\
	\item[(ii)] $-h(-r) \leq |f(z)| \leq h(r)$ \\
	\item[(iii)] $f(\mathbb{D})\supseteq \{w:|w|\leq -h(-1)\}$.
\end{enumerate}
The results are sharp.
\end{lem}

\begin{lem} \cite{Ravichandran-2004} \label{him-p3-lem-1.14}
	Let $\min _{|z|=r} |\phi (z)|=\phi (-r)$, $\max _{|z|=r}|\phi (z)|=\phi (r)$, $|z|=r$. If $f \in \mathcal{C}_{c}(\phi)$, then 
	\begin{enumerate}
		\item [(i)] $k'(-r) \leq |f'(z)| \leq k'(r)$ \\
		\item[(ii)] $-k(-r) \leq |f(z)| \leq k(r)$ \\
		\item[(iii)] $f(\mathbb{D})\supseteq \{w:|w|\leq -k(-1)\}$.
	\end{enumerate}
	The results are sharp.
\end{lem}

Motivated by the class $\mathcal{S}^{*}_{s}$, Gao and Zhou \cite{gao-2005} have studied the class $\mathcal{K}_{s}$ of close-to-convex univalent functions, where $\mathcal{K}_{s}$ is the class of functions $f \in \mathcal{S}$ satisfying the condition 
$$
\real \left(\dfrac{z^{2}f'(z)}{g(z)g(-z)}\right)<0, \quad z\in \mathbb{D}.
$$
A more general class $\mathcal{K}_{s}(\phi)$ has been studied extensively by Cho {\it et.al.} \cite{cho-2011} and Wang {\it et.al.} \cite{wang-2006-a}. For the brevity, we write the definition. 
\begin{defn} \cite{wang-2006-a} \label{him-p3-def-1.2}
For a function $\phi$ with positive real part, the class  $\mathcal{K}_{s}(\phi)$ consists of functions $f \in \mathcal{A}$ satisfying 
$$
-\frac{z^{2}f'(z)}{g(z)g(-z)} \prec \phi (z) \quad \mbox{in} \quad \mathbb{D}
$$
for some function $g \in \mathcal{S}^{*}(1/2)$.
\end{defn} 
In particular, for $\phi(z)=(1+(1-2 \gamma)z)/(1-z)$ with $0\leq \gamma<1$, the class  $\mathcal{K}_{s}(\phi)$ reduces to $\mathcal{K}_{s}(\gamma)$ which has recently been investigated by Kowalczyk and Les-Bomba \cite{kowalczyk-2010}. When $\gamma=0$, we can obtain $\mathcal{K}_{s}$, the subclass of close-to-convex functions which has been defined by Gao and Zhou \cite{gao-2005}. When $\phi (z)=(1+\beta z)/(1-\alpha \beta z)$, where $0 \leq \alpha \leq 1$ and $0< \beta \leq 1$, the class $\mathcal{K}_{s}(\phi)$ reduces to  $\mathcal{K}_{s}(\alpha, \beta)$ defined in \cite{wang-2006-a}.
Now let $q(z)=\sum_{n=1}^{\infty} q_{n}z^{n}$ be analytic in $\mathbb{D}$. Then for fixed $f \in \mathcal{K}_{s}(\phi)$, we define
\begin{equation} \label{him-p3-e-1.20}
S_{f}^{\mathcal{K}}(\phi):= \left\{q(z)=\sum_{n=1}^{\infty} q_{n}z^{n}:\, q \prec\, f \right\}.
\end{equation}
The distortion and growth theorems for the class $\mathcal{K}_{s}(\phi)$ have been obtained in \cite{cho-2011}.
\vspace{1mm}

Let $\phi$ be a Ma-Minda function.
\begin{lem} \cite{cho-2011} \label{him-p3-lem-1.15}
Let $\min _{|z|=r} |\phi (z)|=\phi (-r)$, $\max _{|z|=r}|\phi (z)|=\phi (r)$, $|z|=r$. If $f \in \mathcal{K}_{s}(\phi)$, then the following sharp inequalities hold:
\begin{enumerate}
	\item [(i)] $\dfrac{\phi (-r)}{1+r^{2}} \leq |f'(z)| \leq \dfrac{\phi (r)}{1-r^{2}}$    
	   $(|z|=r<1)$ \\[3mm]
	\item [(ii)] $ \displaystyle{\int \limits _{0}^{r}}\dfrac{\phi (-t)}{1+t^{2}} \,\, dt\, \leq |f(z)| \leq \displaystyle{\int \limits _{0}^{r}}\dfrac{\phi (t)}{1-t^{2}}\,\, dt$ $(\,|z|=r<1)$.
\end{enumerate}
\end{lem}

\vspace{4mm}

 Let $f$ and $g$ be two analytic functions in $\mathbb{D}$ such that $g \prec f$. Let 
\begin{equation} \label{him-p2-e-1.22}
g(z)=\sum_{n=0}^{\infty} b_{n}z^{n}.
\end{equation}
In 2018, Bhowmik and Das \cite{bhowmik-2018} proved the following interesting result for subordination classes.
\begin{lem} \label{him-p2-lem-1.23} \cite{bhowmik-2018}
	Let $f$ and $g$ be analytic in $\mathbb{D}$ with Taylor expansions \eqref{him-p3-e-1.1} and \eqref{him-p2-e-1.22} respectively and $g \prec f$, then 
	\begin{equation} \label{him-p2-e-1.24}
	\sum_{n=0}^{\infty} |b_{n}| r^{n} \leq \sum_{n=0}^{\infty} |a_{n}| r^{n}
	\end{equation}
	for $z|=r \leq 1/3.$
\end{lem} 
\vspace{6mm}
In general, one obtains the Bohr radius for certain classes of analytic functions in $\mathbb{D}$, when the sharp coefficient bounds for this class are known.
% For example, consider $\mathcal{S}^{*}(\alpha)$, the class of starlike functions of order $\alpha$ and $\mathcal{C}(\alpha)$, the class of convex functions of order $\alpha$ and so on.
But the sharp coefficient bounds for most of the Ma-Minda subclasses are not yet known.
Using Lemma \ref{him-p2-lem-1.23}, Allu and Halder \cite{Himadri-Vasu-P2} recently have obtained Bohr radius for certain classes of Ma-Minda starlike and convex functions. In this article, we consider certain classes of close-to-convex functions associated with Ma-Minda functions {\it e.g.} $\mathcal{S}_{c}^{*}(\phi), \, \mathcal{C}_{c}(\phi),\, \mathcal{K}_{s}(\phi)\quad \mbox{and} \quad \mathcal{C}_{s}(\phi)$. The sharp coefficient bounds of these classes are not yet known. Hence, we encounter the problem to find the best possible lower bound of the radius so that Bohr phenomenon holds for these classes. As a consequence, we also establish the Bohr phenomenon for several important subclasses for particular choices of $\phi$.
\section{Main Results}
Before going to state our main results we prove an preliminary result which is required to prove some of our results.

\begin{lem} \label{him-p3-lem-2.1}
\begin{enumerate}
\item[(i)] Let  $f$ and $g$ be analytic in $\mathbb{D}$ with series representation $f(z)=\sum_{n=1}^{\infty} a_{n}z^{n}$ and \eqref{him-p2-e-1.22} respectively such that $f(z)=\int_{0}^{z}g(\xi)\, d\xi$ for $z \in \mathbb{D}$, where integration is taken along a linear segment joining $0$ to $z \in \mathbb{D}$. Then 
$$
M_{f}(r) = \int_{0}^{r}M_{g}(t) \, dt \quad \mbox{for } \quad |z|=r<1.
$$
Here $M_{f}(r)$ and $M_{g}(r)$ are respectively the majorant series associated with $f$ and $g$ respectively.
\item[(ii)] 	Let $f$ and $g$ be analytic in $\mathbb{D}$ with Taylor expansions \eqref{him-p3-e-1.1} and \eqref{him-p2-e-1.22} respectively and $g \prec f$, then 
$M_{G}(r) \leq M_{F}(r)$ for $|z|=r\leq 1/3$, where $G(z)=\int_{0}^{z}g(\xi)\, d\xi$ and $F(z)=\int_{0}^{z}f(\xi)\, d\xi$ for $z \in \mathbb{D}$.
\end{enumerate}
\end{lem}

Let $\min _{|z|=r} |\phi (z)|=\phi (-r)$ and $\max _{|z|=r}|\phi (z)|=\phi (r)$, $|z|=r$. We assume these notations throught this paper. Here $\phi$ is the Ma-Minda function.
\begin{thm} \label{him-p3-thm-2.2}
Let $f \in \mathcal{K}_{s}(\phi)$ be of the form \eqref{him-p3-e-1.5}. Then 
\begin{equation} \label{him-p3-e-2.2-a}
|z|+ \sum_{n=2}^{\infty} |a_{n}| |z|^{n} \leq d(f(0),\partial f(\mathbb{D})) 
\end{equation}
for $|z|=r \leq R_{f}$, where $R_{f}=\min \{1/3,r_{f}\}$ and $r_{f}$ is the smallest positive root of $R(r)=L(1)$ in $(0,1)$. Here $R(r):=\int  _{0}^{r}\left(M_{\phi} (t))/(1-t^{2}\right) \, dt$ , $L(r):=\int _{0}^{r}\left(\phi (-t)\right)/(1+t^{2}) \, dt$ and $M_{\phi}$ is the associated majorant series of $\phi$.
\end{thm}
\begin{rem} \label{him-p3-rem-2.1}
\begin{enumerate}
\item [(i)] Assume that the coeficients of $\phi(z)=1+\sum_{n=1}^{\infty} B_{n}z^{n}$ in the Theorem \ref{him-p3-thm-2.2} are all positive {\it i.e.} $B_{n}>0$ for $n\geq 1$. Then the majorant series $M_{\phi}(r)=\phi (r)$, $0<r<1$ and hence $R(r):=\int  _{0}^{r}(\phi (t))/(1-t^{2}) \, dt$.
% When $r_{f}\leq 1/3$, the constant $R_{f}$ is the best possible. 
\item [(ii)] (Bohr phenomenon for the corresponding class $\mathcal{K}_{s}(\Phi)$ associated with non-Ma-Minda functions) 
 Let $\Phi$ be the corresponding non-Ma-Minda function of $\phi$, which is actually a rotation by mere replacing $z$ by $-z$. Therefore the image of the unit disk $\mathbb{D}$ under the functions $\Phi$ and $\phi$ are identical. Thus we conclude that $\mathcal{K}_{s}(\Phi)=\mathcal{K}_{s}(\phi)$ and the Bohr phenomenon \eqref{him-p3-e-2.2-a} holds for the class $\mathcal{K}_{s}(\Phi)$ for the same $R_{f}$.
\end{enumerate}
\end{rem}

{\bf Some applications:}
\begin{lem} \label{him-p3-lem-2.5} $( \mbox{Bohr phenomenon for the corresponding subordination class} )$ \\
	Let $q(z) =\sum  _{n=1}^{\infty} q_{n}z^{n}\in S_{f}^{\mathcal{K}}(\phi)$ as defined in \eqref{him-p3-e-1.20} and $f$ be of the form \eqref{him-p3-e-1.5}. Then 
	$$
	\sum_{n=1}^{\infty} |q_{n}| |z|^{n} \leq d(f(0),\partial f(\mathbb{D})) 
	$$ 
	for $|z|=r \leq R_{f}$, where $R_{f}$ is defined as in Theorem \ref{him-p3-thm-2.2}.
\end{lem}
For $\phi (z)=(1+(1-2\gamma)z)/(1-z)$, the class $\mathcal{K}_{s}(\phi)$ reduces to $\mathcal{K}_{s}(\gamma)$. In particular, for $\gamma=0$, $\mathcal{K}_{s}(\phi)$ reduces to $\mathcal{K}_{s}$.
\begin{cor} \label{him-p3-cor-2.7}
\begin{enumerate} 
	\item [(i)] $(\mbox{Bohr phenomenon for the class $\mathcal{K}_{s}(\gamma)$})$\\
	 Any function $f \in \mathcal{K}_{s}(\gamma)$ with $0 \leq \gamma <0.259056404$ satisfies the inequality \eqref{him-p3-e-2.2-a} for $|z|=r \leq r_{f}$, where $r_{f}$ is the root of 
	\begin{equation} \label{him-p3-e-2.6}
	\frac{\gamma}{2}ln\left(\frac{1+r}{1-r}\right)+(1-\gamma)\frac{r}{1-r}=\frac{1-\gamma}{2}ln2+\frac{\gamma \pi}{4}
	\end{equation}
	in $(0,1/3)$.
	
%	The radius $r_{f}$ is the best possible.
	\item [(ii)] Each function $f \in \mathcal{K}_{s}$ satisfies the Bohr inequality \eqref{him-p3-e-2.2-a} for $|z|=r \leq r_{f}$, where $r_{f}=ln \, 2/(2+ln \, 2) \approx 0.257374415$. 
	
\end{enumerate}
\end{cor}
For $\phi (z)=(1+\beta z)/(1-\alpha \beta z)$, where $0 \leq \alpha \leq 1$ and $0< \beta \leq 1$, the class $\mathcal{K}_{s}(\phi)$ reduces to  $\mathcal{K}_{s}(\alpha, \beta)$. In particular, for $\alpha=\beta=1$, $\mathcal{K}_{s}(\alpha, \beta)$ coincides with the class $\mathcal{K}_{s}$.
\begin{cor} \label{him-p3-cor-2.8}
 The class $\mathcal{K}_{s}(\alpha, \beta)$ satisfies the Bohr phenomenon \eqref{him-p3-e-2.2-a} for $|z|=r \leq R_{f}= \min \{1/3,r_{f}\}$, where $r_{f}$ is the smallest root of 
	\begin{equation} \label{him-p3-e-2.10}
	\int \limits _{0}^{r} \dfrac{1+\beta t}{(1-\alpha \beta t)(1-t^{2})}\, dt = \int \limits _{0}^{1} \dfrac{1-\beta t}{(1+\alpha \beta t)(1+t^{2})}\, dt 
	\end{equation}
	in $(0,1)$.
	%  and we obtain $r_{f}$ by putting $\alpha=\beta=1$ in \eqref{him-p3-e-2.10}, is the smallest root of 
%	$$
%	\int \limits _{0}^{r} \dfrac{1+t}{(1- t)(1-t^{2})}\, dt = \int \limits _{0}^{1} \dfrac{1-t}{(1+ t)(1+t^{2})}\, dt 
%	$$
%	in $(0,1)$.
\end{cor}

\begin{thm} \label{him-p3-thm-2.11}
Let $f \in \mathcal{S}_{c}^{*}(\phi)$ be of the form \eqref{him-p3-e-1.5}. Then 
\begin{equation} \label{him-p3-e-2.11-a}
|z|+ \sum_{n=2}^{\infty} |a_{n}| |z|^{n} \leq d(f(0),\partial f(\mathbb{D}))
\end{equation}
for $|z|=r \leq \min \{1/3,r_{f}\}$ and $r_{f}$ is the smallest positive root of 
$P(r)+h(-1)=0$ in $(0,1)$, where $P(r):=\int  _{0}^{r}\left(\left(M_{h}(t) M_{\phi}(t)\right)/t \right)\, dt$. Here $M_{h}(t)$ and $M_{\phi}(t)$ are the majorant series of $h$ and $\phi$ respectively.
\end{thm}

\begin{rem} \label{him-p3-rem-2.2}
\begin{enumerate}
\item [(i)] $\left(\mbox{Bohr radius for $\mathcal{S}_{c}^{*}(\phi)$ when $\phi$ has positive coefficients}\right)$ \\
 Let $\phi(z)=1+\sum_{n=1}^{\infty} B_{n}z^{n}$. It is worth to point out that if we impose an additional condition on $\phi$ that the coefficients $B_{n}$'s are positive, then the majorant series $M_{\phi}(r)=\phi(r)$. From the definition of $h$ in \eqref{him-p3-e-1.6}, we obtain 
 \begin{equation} \label{him-p2-e-3.8}
 h(z)= z \exp \left(\int\limits_{0}^{z} \dfrac{\phi(t)-1}{t}\,\, dt\right)=z \exp \left(\sum \limits _{n=1}^{\infty} \frac{B_{n}}{n}z^{n}\right).
 \end{equation}
 From \eqref{him-p2-e-3.8}, it is easy to see that
  $$M_{h}(r)=h(r) \quad \mbox{and} \quad 
P(r)=\int_{0}^{r} \left(\left( h(t)\phi(t)\right)/t \right) \, dt =h(r).
$$
Then each $f \in \mathcal{S}_{c}^{*}(\phi)$ satisfies the inequality \eqref{him-p3-e-2.11-a} for $|z|\leq \min \{1/3,r_{f}\}$, where $r_{f}$ is the root of the equation $h(r)+h(-1)=0$. In particular, when $r_{f} \leq 1/3$, the radius $r_{f}$ is the best possible for the function $f=h \in \mathcal{S}_{c}^{*}(\phi)$, since it has real coefficients and belongs to $\mathcal{S}^{*}(\phi)$. Indeed, for $|z|=r_{f}$, $M_{h}(r_{f})=h_{r_{f}}=-h(-1)=d(h(0),\partial h(\mathbb{D}))$, which shows that $r_{f}$ is the best possible.
\item [(ii)] (Bohr phenomenon for corresponding class $\mathcal{S}^{*}_{c}(\Phi)$ associated with non-Ma-Minda function) 
Let $\Phi$ be the corresponding non-Ma-Minda function of $\phi$. Since $\Phi$ is actually obtained from $\phi$ by a rotation $z$ by $-z$, the image of the unit disk $\mathbb{D}$ under the functions $\Phi$ and $\phi$ are identical. Thus we conclude that $\mathcal{S}^{*}_{c}(\Phi)=\mathcal{S}^{*}_{c}(\phi)$ and the Bohr radius for the class $\mathcal{S}^{*}_{c}(\Phi)$ is same as that of $\mathcal{S}^{*}_{c}(\phi)$.
\end{enumerate}
\end{rem}

Let $S^{*}_{cf}(\phi)$ denote the class of analytic functions $g$ which are subordinate to a fixed function $f \in \mathcal{S}^{*}_{c}(\phi)$.
\begin{lem} \label{him-p3-lem-2.13} $\left(\mbox{Bohr phenomenon for the corresponding subordination class $S^{*}_{cf}(\phi)$} \right)$
Let $g \in S^{*}_{cf}(\phi)$ be of the form $g(z)=\sum_{n=1}^{\infty} g_{n}z^{n}$. Then 
\begin{equation} \label{him-p3-e-2.13-a}
\sum_{n=1}^{\infty} |g_{n}| |z|^{n} \leq d(f(0),\partial f(\mathbb{D}))
\end{equation}
for $|z|=r \leq \min \{1/3,r_{f}\}$, where $r_{f}$ is given as in  Theorem \ref{him-p3-thm-2.11}.
\end{lem}
Similar results on Bohr phenomenon for the class $\mathcal{S}^{*}_{c}(\phi)$ hold for the class $S^{*}_{cf}(\phi)$. In view of the Remark \ref{him-p3-rem-2.2} and Lemma \ref{him-p3-lem-2.13}, we obtain the following interesting corollaries. Let $\phi (z)= (1+sz)^{2}$ with $0.444981<s \leq 1/\sqrt{2}$, then $\mathcal{S}^{*}_{c}(\phi)$ reduces to the class $\mathcal{S}_{c}^{*} \left((1+sz)^{2}\right)$.
\begin{cor} \label{him-p3-cor-2.15}
 The class $ \mathcal{S}_{c}^{*} \left((1+sz)^{2}\right)\left(\mbox{and} \quad \mathcal{S}_{cf}^{*} \left((1+sz)^{2}\right)\right)$ satisfies the Bohr inequality  \eqref{him-p3-e-2.11-a} for $|z|=r \leq r_{f}$, where  $0<r_f<1/3$ and $r_{f}$ is the root of the equation 
\begin{equation} \label{him-p3-e-2.16}
r\exp \left(s\left(2r+\dfrac{sr^{2}}{2}\right)\right)=\exp \left(s\left(-2+\dfrac{s}{2}\right)\right).
\end{equation}
 The radius $r_{f}$ is the best posible.
\end{cor}
\begin{table}[ht]
	
	\begin{tabular}{|l|l|} 
		\hline
		$s$& $r_{f}$ \\
		\hline
		$0.1$& $0.71184$\\
		\hline
		$0.15$& $0.619461$\\
		\hline
		$0.2$ & $0.546344$\\
		\hline
		$0.25$& $0.486934$\\
		\hline
		$0.3$& $0.437693$\\
		\hline
		$0.35$& $0.39624$\\
		\hline
		$0.4$& $0.360903$\\
		\hline
	\end{tabular}
\hspace{16mm}
\begin{tabular}{|l|l|} 
	\hline
	$s$& $r_{f}$ \\
	\hline
	$0.45$& $0.330472$\\
	\hline
	$0.5$& $0.3040402$\\
	\hline
	$0.55$& $0.28091732$\\
	\hline
	$0.6$& $0.2605657$\\
	\hline
	$0.65$& $0.24256$\\
	\hline
	$0.7$& $0.226558$\\
	\hline
	$1/\sqrt{2}$& $0.22443096$\\
	\hline
\end{tabular}
\vspace{3mm}
\caption{The radius $r_{f}$ for different values of $s$}
\label{tabel-1}
\end{table}
From Table \ref{tabel-1}, it is easy to see that $r_{f}>1/3$ when $s<0.444981$ and hence Bohr phenomenon holds for $r \leq 1/3$ and  $r_{f}<1/3$ when $0.444981<s\leq 1/\sqrt{2}$. Therefore the radius $r_{f}$ is the best possible.
\begin{cor} \label{him-p3-cor-2.17}
For $\phi(z)=\alpha +(1-\alpha) e^{z}$ with $0 \leq \alpha <0.05284$, the class $\mathcal{S}^{*}_{c}(\phi)$ satisfies the Bohr phenomenon \eqref{him-p3-e-2.11-a} for $|z|=r \leq r_{f}$, where  $0<r_f<1/3$ . The rdius $r_{f}$ is the best possible. 

\begin{table}[ht]
	\centering
	\begin{tabular}{|l|l|l|l|l|}
		\hline
		$\alpha$& $h(1/3)$& $h(-1)$& Sign of $D_{2}(0)$ & Sign of $D_{2}(1/3)$ \\
		\hline
		$0.0$& $0.47935$& $0.4508594$& $-$& $+$\\
		\hline
		$0.01$& $0.477619$& $0.454465$& $-$& $+$\\
		\hline
		$0.02$& $0.475887697$& $0.458100015$& $-$& $+$\\
		\hline
		$0.03$& $0.47416191$& $0.4617638$& $-$& $+$\\
		\hline
		$0.04$& $0.47244238$& $0.465456$& $-$& $+$\\
		\hline
		$0.05$& $0.470729$& $0.469179$& $-$& $+$\\
		\hline
		$0.06$& $0.469022$& $0.47293$& $-$& $-$\\
		\hline
		$0.07$& $0.46732112$& $0.4767143$& $-$& $-$\\
		\hline
	\end{tabular}
    \vspace{3mm}
    \caption{Existance of the sharp radius $r_{f}$ in $(0,1/3)$ for different values of $\alpha$ in $[0,0.05284)$}
    \label{tabel-2}
\end{table}
\end{cor}
From Table \ref{tabel-2}, it is clear that $r_{f}$ lies in $(0,1/3)$ when $0 \leq \alpha <0.05284$ and hence $r_{f}$ is the best posiible. On the other hand $r_{f}>1/3$ for $\alpha > 0.05284$ and the corresponding Bohr phenomenon holds for $r\leq 1/3$.
\begin{cor} \label{him-p3-cor-2.18}
Let $\phi(z)=\left((1+z)/(1-z)\right)^{\alpha}$ with $0< \alpha \leq 1$. Also assume $h(1/3)>-h(-1)$, where 
$$
h(r)=r\exp \left(\int \limits _{0}^{r} \frac{\left(\frac{1+t}{1-t}\right)^{\alpha}-1}{t}\,\, dt \right)
$$
and 
$$
-h(-1)= \exp \left(\int \limits _{0}^{-1} \frac{\left(\frac{1+t}{1-t}\right)^{\alpha}-1}{t}\,\, dt \right).
$$
Then the class $\mathcal{S}^{*}_{c}(\phi)$ satisfies the Bohr phenomenon \eqref{him-p3-e-2.11-a} for $|z|=r \leq r_{f}$, where $r_{f}$ is the smallest root of the equation $D_{3}(r):=h(r)+h(-1)=0$.
\end{cor}

\begin{table}[ht]
	\centering
	\begin{tabular}{|l|l|l|l|l|}
		\hline
		$\alpha$& $h(1/3)$& $-h(-1)$& Sign of $D_{3}(0)$ & Sign of $D_{3}(1/3)$ \\
		\hline
		$0.2$& $0.38335$& $0.65515$& $-$& $-$\\
		\hline
		$0.4$& $0.4453711$& $0.475453$& $-$& $-$\\
		\hline
		$0.45$& $0.4631699$& $0.443795$& $-$& $+$\\
		\hline
		$0.5$& $0.482023$& $0.415759$& $-$& $+$\\
		\hline
		$0.6$& $0.523214$& $0.368431$& $-$& $+$\\
		\hline
		$0.7$& $0.569663$& $0.330139$& $-$& $+$\\
		\hline
		$0.8$& $0.62222$& $0.298621$& $-$& $+$\\
		\hline
		$0.9$& $0.681928$& $0.272286$& $-$& $+$\\
		\hline
	\end{tabular}
    \vspace{3mm}
    \caption{Existance of the sharp radius $r_{f}$ in $(0,1/3)$ for different values of $\alpha$}
    \label{table-3}
\end{table}
From the Table \ref{table-3}, it is evident that for different values of $\alpha$, the constant $r_{f}$ sometimes does not lie in $(0,1/3)$. However, when $r_{f}$ lies in $(0,1/3)$, the corresponding $r_{f}$ is the best possible and the Bohr phenomenon for the class $\mathcal{S}^{*}_{c}(\phi)$ holds for $r\leq r_{f}$.
\begin{cor} \label{him-p3-cor-2.19}
Let $\phi (z)=\left(1+(1-2\gamma)z\right)/(1-z)$ with $0 \leq \gamma <1/2$. Then each $f \in \mathcal{S}_{c}^{*}\left(\left(1+(1-2\gamma)z\right)/(1-z) \right)$ satisfies the inequality \eqref{him-p3-e-2.11-a} for $|z|=r \leq r_{f}$, where  $0<r_{f}<1/3$ and $r_{f}$ is the root of 
\begin{equation} \label{him-p3-e-2.20}
r+2r^{1/(2(1-\gamma))}-1=0.
\end{equation}
The radius $r_{f}$ is the best possible.
\end{cor}

\begin{cor} \label{him-p3-cor-2.21}
If $\phi (z)=(1+Az)/(1+Bz)$ with $-1 \leq B<A\leq 1$, then 
\begin{enumerate}
	\item [(i)] when $B=0$, every function $f \in \mathcal{S}_{c}^{*}\left((1+Az)/(1+Bz)\right)$ satisfies the inequality \eqref{him-p3-e-2.11-a} for $|z|=r \leq r_{f}$, where  $0<r_{f}<1/3$ and $r_{f}$ is the unique root of
	\begin{equation} \label{him-p3-e-2.22}
	re^{Ar}=e^{-A},
	\end{equation}
	provided $A\geq (3/4)ln\, 3$.
	The radius $r_{f}$ is the best possible.
	\item [(ii)] When $B\neq 0$, every function $f \in \mathcal{S}_{c}^{*}\left((1+Az)/(1+Bz)\right)$ satisfies the inequality \eqref{him-p3-e-2.11-a} for $|z|=r \leq r_{f}$, where  $0<r_{f}<1/3$ and $r_{f}$ is the unique root of
	\begin{equation} \label{him-p3-e-2.23}
	r\left(1+Br\right)^{\frac{A-B}{B}}=\left(1-B\right)^{\frac{A-B}{B}},
	\end{equation}
	provided $(1/3) \left(1+B/3\right)^{(A-B)/B} \geq \left(1-B\right)^{(A-B)/B}$.
	The radius $r_{f}$ is the best possible.
\end{enumerate}
\end{cor}

\begin{table}[ht]
	\centering
	$(A=1)$
	\begin{tabular}{|l|l|}
		\hline
		$B$& $r_{f}$ \\
		\hline
		$-0.1$& $0.261789$\\
		\hline
		$-0.2$& $0.247088$\\
		\hline
		$-0.3$ & $0.23402$\\
		\hline
		$-0.4$& $0.222323$\\
		\hline
		$-0.5$& $0.21179$\\
		\hline
		$-0.6$& $0.202239$\\
		\hline
		$-0.7$& $0.193548$\\
		\hline
		$-0.8$& $0.185599$\\
		\hline
		$-0.9$& $0.1783$\\
		\hline
		$-1.0$& $0.17157$\\
		\hline
	\end{tabular}
\hspace{10mm} $(A=1/2)$
\begin{tabular}{|l|l|}
	\hline
	$B$& $r_{f}$ \\
	\hline
	$-0.1$& $0.432852$\\
	\hline
	$-0.2$& $0.395824$\\
	\hline
	$-0.3$ & $0.364714$\\
	\hline
	$-0.4$& $0.338205$\\
	\hline
	$-0.5$& $0.31534$\\
	\hline
	$-0.6$& $0.295418$\\
	\hline
	$-0.7$& $0.277899$\\
	\hline
	$-0.8$& $0.262372$\\
	\hline
	$-0.9$& $0.248514$\\
	\hline
	$-1.0$& $0.236068$\\
	\hline
\end{tabular}
\vspace{3mm}
\caption{The radius $r_{f}$ for different values of $B$ when $A=1$ and $A=1/2$}
\end{table}
From the Table $4$, we see that for different values of $A$ and $B$, sometimes the radius $r_{f}<1/3=0.33333$ and in that case $r_{f}$ is the best possible. When $r_{f}>1/3$, Bohr phenomenon for class $\mathcal{S}_{c}^{*}\left((1+Az)/(1+Bz)\right)$ holds for $r \leq 1/3$.

\begin{thm} \label{him-p3-thm-2.24}
	Let $f \in \mathcal{C}_{c}(\phi)$ be of the form \eqref{him-p3-e-1.5}. Then 
	\begin{equation} \label{him-p3-e-2.25}
	|z|+ \sum_{n=2}^{\infty} |a_{n}| |z|^{n} \leq d(f(0),\partial f(\mathbb{D}))
	\end{equation}
	for $|z|=r \leq \min \{1/3,r_{f}\}$ and $r_{f}$ is the smallest positive root of 
	$T(r)=-k(-1)$ in $(0,1)$ and $$T(r):=\int \limits _{0}^{r}\frac{1}{s}\int \limits _{0}^{s}M_{k'}(t) M_{\phi}(t) \, dt\, ds.$$ 
	Here $M_{k'}(t)$ and $M_{\phi}(t)
	$ are respectively the majorant series of $k'$ and $\phi$ respectively.
\end{thm}

\begin{thm} \label{him-p3-thm-2.25}
Let $f \in \mathcal{C}_{s}(\phi)$ be of the form \eqref{him-p3-e-1.5}. Then 
\begin{equation} \label{him-p3-e-2.25-a}
|z|+ \sum_{n=2}^{\infty} |a_{n}| |z|^{n} \leq d(f(0),\partial f(\mathbb{D}))
\end{equation}
for $|z|=r \leq \min \{1/3,r_{f}\}$ and $r_{f}$ is the smallest positive root of 
$R_{s}(r)=L_{s}(1)$ in $(0,1)$, where 
$$
R_{s}(r):=\int \limits _{0}^{r} \frac{1}{s} \int \limits _{0}^{s} M_{K'}(t) M_{\phi}(t) \, dt\, ds \quad \mbox{and} \quad L_{s}(r):=\int \limits _{0}^{r} \frac{1}{s} \int \limits _{0}^{s} \left(k'(-t^{2})\right)^{1/2}\phi(-t) \, dt\, ds
$$
and $K'(r)=\left(k'(t^{2})\right)^{1/2}$. 
\end{thm}
\begin{rem}
	\begin{enumerate}
		\item [(i)] Let $\Phi$ be the corresponding non-Ma-Minda class with respect to $\phi$. Then the Bohr radius for the class $\mathcal{C}_{s}(\Phi)$ is same as that of $\mathcal{C}_{s}(\phi)$.
		\item [(ii)] Let $S^{*}_{sf}(\phi)$ be the class of analytic functions $g$ of the form $g(z)=\sum_{n=1}^{\infty} g_{n}z^{n}$ in $\mathbb{D}$ subordinate to a fixed function $f \in \mathcal{C}_{s}(\Phi)$, then 
		$$
		\sum \limits _{n=1}^{\infty} |g_{n}| |z|^{n} \leq d(f(0),\partial f(\mathbb{D}))
		$$ 
		for $|z|=r \leq \min \{1/3,r_{f}\}$ and $r_{f}$ is given as in Theorem \ref{him-p3-thm-2.25}.
	\end{enumerate}
\end{rem}

\section{Proof of the main results}
\begin{pf} [{\bf Proof of Lemma   \ref{him-p3-lem-2.1}}]
\begin{enumerate}
	\item [(i)]  In view of the relation $f(z)=\int_{0}^{z}g(\xi)\, d\xi$, we obtain 
	$$
	\sum_{n=1}^{\infty} a_{n}z^{n}=\sum_{n=1}^{\infty} \frac{b_{n-1}}{n}z^{n}.
	$$
	Therefore 
	$$M_{f}(r)=\sum \limits _{n=1}^{\infty} \frac{|b_{n-1}|}{n}r^{n}= \int \limits  _{0}^{r} \sum_{n=0}^{\infty} |b_{n}|t^{n} \, dt =\int_{0}^{r}M_{g}(t) \, dt \quad \mbox{for} \quad r<1.$$
	\item[(ii)] From Lemma \ref{him-p2-lem-1.23}, we have $M_{g}(r) \leq M_{f}(r)$ for $r \leq 1/3$ and integrating this we obtain 
	$$
	\int \limits _{0}^{r}M_{g}(t) \, dt \leq \int \limits _{0}^{r}M_{f}(t) \, dt \quad \mbox{for} \quad r\leq 1/3.
	$$
	Hence from the first part of this Lemma, we obtain 
	$$
	M_{G}(r)=\int \limits _{0}^{r}M_{g}(t) \, dt \leq \int \limits _{0}^{r}M_{f}(t) \, dt= M_{F}(r)
	$$	
\end{enumerate}
for $r \leq 1/3$.
\end{pf}

\begin{pf} [{\bf Proof of Theorem   \ref{him-p3-thm-2.2}}]
Let $f \in \mathcal{K}_{s}(\phi)$, then from Lemma \ref{him-p3-lem-1.15}, the Euclidean distance between $f(0)$ and the boundary of $f(\mathbb{D})$ is  
\begin{equation} \label{him-p3-e-3.1-a}
d(f(0), \partial f(\mathbb{D}))= \liminf \limits_{|z|\rightarrow 1} |f(z)-f(0)| \geq \int \limits _{0}^{1}\dfrac{\phi (-t)}{1+t^{2}} \,\, dt.
\end{equation}
 By the subordination principle, there exists an analytic function $\omega :\mathbb{D} \rightarrow \mathbb{D}$ with $\omega(0)=0$ such that
\begin{equation} \label{him-p3-e-3.1}
-\dfrac{z^{2}f'(z)}{g(z)g(-z)}=\phi(\omega(z)).
\end{equation}
Let $G(z):=-g(z)g(-z)/z$. Clearly, $G$ is an odd starlike function in $\mathbb{D}$.  Let $G(z)=z+\sum_{n=2}^{\infty}g_{2n-1}z^{2n-1}$. It is well-known that $|g_{2n-1}|\leq 1$ for $n \geq 2$. Therefore 
\begin{equation} \label{him-p3-e-3.2}
M_{G}(r) \leq r+\sum_{n=2}^{\infty}r^{2n-1}=\frac{r}{1-r^{2}}, \quad 0<r<1.
\end{equation}
From \eqref{him-p3-e-3.1}, we have $zf'(z)=G(z)\phi(\omega(z))$, which immediately follows that 
\begin{equation} \label{him-p3-e-3.2-a}
f(z)=\int \limits_{0}^{z} \dfrac{G(\xi)\phi(\omega (\xi))}{\xi} \,\, d\xi.
\end{equation}
It is known that for two analytic functions $f$ and $g$ in $\mathbb{D}$, $M_{fg}(r)\leq M_{f}(r)M_{g}(r)$, where $M_{f}(r)$, $M_{g}(r)$ and $M_{fg}(r)$ are associated majorant series of $f$, $g$ and the product $fg$ respectively. Therefore $M_{G(\phi \circ \omega)}(r)\leq M_{G}(r)M_{\phi \circ \omega}(r)$. Since $\phi \circ \omega \prec \phi$, by Lemma \ref{him-p2-lem-1.23}, we have 
\begin{equation} \label{him-p3-e-3.2-b}
M_{\phi \circ \omega} (r) \leq M_{\phi} (r) \quad \mbox{for}\quad |z|=r\leq 1/3.
\end{equation}
In view of Lemma \ref{him-p3-lem-2.1} and \eqref{him-p3-e-3.2}, \eqref{him-p3-e-3.2-a} and \eqref{him-p3-e-3.2-b}, we obtain
\begin{equation} \label{him-p3-e-3.5}
M_{f}(r) \leq \int \limits_{0}^{r} \dfrac{M_{G}(t)M_{\phi \circ \omega}(t)}{t} \, \, dt \leq \int \limits _{0}^{r}\frac{M_{\phi} (t)}{1-t^{2}} \,\, dt=R(r)
\end{equation}
for $|z|=r\leq 1/3$. We note that $R(r) \leq L(1)$ whenever $r \leq r_{f}$, where $r_{f}$ is the smallest positive root of $R(r)=L(1)$ in $(0,1)$. Let $H_{1}(r)=R(r)-L(1)$ then $H_{1}(r)$ is continuous function in $[0,1]$. Since $R(1)>L(1)$ and $M_{\phi}(t) \geq |\phi(t)|$, it follows that   
$$
H_{1}(0)=L(1)=-\int \limits _{0}^{1}\dfrac{\phi (-t)}{1+t^{2}} \,\, dt <0
$$
and 
$$
H_{1}(1)=R(1)-L(1)=\int \limits _{0}^{1}\dfrac{M_{\phi} (t)}{1-t^{2}} \, \, dt -\, \int \limits _{0}^{1}\dfrac{\phi (-t)}{1+t^{2}} \,\, dt\, >0.
$$
 Therefore $H_{1}$ has a root in $(0,1)$. Let $r_{f}$ be the smallest root of $H_{1}$ in $(0,1)$. Then $R(r)\leq L(1)$ for $r\leq r_{f}$.
%Therefore, we obtain 
%\begin{equation} \label{him-p3-e-3.6}
%R(r_{f})=L(1)=\int \limits _{0}^{1}\dfrac{\phi (-t)}{1+t^{2}} \, dt.
%\end{equation} 
From \eqref{him-p3-e-3.1-a} and \eqref{him-p3-e-3.5}, we obtain 
$$
M_{f}(r) \leq \int \limits _{0}^{1}\dfrac{\phi (-t)}{1+t^{2}} \,\, dt \leq d(f(0), \partial f(\mathbb{D}))
$$
for $|z|=r \leq \min \{1/3,r_{f}\}=R_{f}$. 
\end{pf}

\begin{pf} [{\bf Proof of Lemma  \ref{him-p3-lem-2.5}}]
	From the definition of $S
	_{f}^{\mathcal{K}}(\phi)$, we have $q \prec f$. In view of Lemma \ref{him-p2-lem-1.23}, we obtain $M_{q}(r) \leq M_{f}(r)$ for $|z|=r \leq 1/3$. Hence from \eqref{him-p3-e-2.2-a}, we get $\sum_{n=1}^{\infty}|q_{n}||z|^{n} \leq d(f(0), \partial f(\mathbb{D}))$ for $|z|=r \leq \min \{1/3, r_{f}\}$.	
\end{pf}
\begin{pf} [{\bf Proof of Corollary   \ref{him-p3-cor-2.7}}]
\begin{enumerate}
\item[(i)] Let $f \in \mathcal{K}_{s}(\gamma)$. Then a simple computation shows that 
$$
R(r)=\frac{\gamma}{2}ln\left(\frac{1+r}{1-r}\right)+(1-\gamma)\frac{r}{1-r}
$$
and 
$$
L(r)=(1-\gamma) ln\left(\frac{1+r}{\sqrt{1+r^{2}}}\right)+ \gamma \arctan r.
$$
Clearly, $L(1)=\left((1-\gamma)/2\right)ln2+ \gamma\pi/4$ and $H_{1}(r):=R(r)-L(1)$. Then $H_{1}$ is continuous in $[0,1)$. A simple computation shows that $H_{1}(0)<0$ and 
$H_{1}(1/3)>0$ if $0 \leq \gamma <0.259056404$. Therefore, $H_{1}$ has a root in $(0,1/3)$ and choose the smallest root to be $r_{f}$ in $(0,1/3)$. Thus the inequality \eqref{him-p3-e-2.2-a} holds for $|z|=r \leq r_{f}$.
\item [(ii)] Putting $\gamma=0$ in \eqref{him-p3-e-2.6}, we obtain $r_{f}=ln\,2/(2+ln\,2)$.
\end{enumerate}
\end{pf}

\begin{pf} [{\bf Proof of Theorem   \ref{him-p3-cor-2.8}}]	
It is easy to see that the coefficients of the power series of $\phi (z)=(1+\beta z)/(1-\alpha \beta z)$ are positive, where $0 \leq \alpha \leq 1$ and $0< \beta \leq 1$. In view of Remark \ref{him-p3-rem-2.1} (i), we obtain $M_{\phi}(r)=\phi (r)$ and 
$$
R(r)=\int \limits _{0}^{r} \dfrac{1+\beta t}{(1-\alpha \beta t)(1-t^{2})}\, \,dt.
$$
Therefore, from Theorem \ref{him-p3-thm-2.2}, $r_{f}$ is the root of 
$$
\int \limits _{0}^{r} \dfrac{1+\beta t}{(1-\alpha \beta t)(1-t^{2})}\,\, dt = \int \limits _{0}^{1} \dfrac{1-\beta t}{(1+\alpha \beta t)(1+t^{2})}\,\, dt .
$$
Thus, the class $\mathcal{K}_{s}(\alpha, \beta)$ satisfies the Bohr phenomenon \eqref{him-p3-e-2.2-a} for $|z|=r \leq R_{f}= \min \{1/3,r_{f}\}$.
\end{pf}
\begin{pf} [{\bf Proof of Theorem   \ref{him-p3-thm-2.11}}]	
Let $f \in \mathcal{S}_{c}^{*}(\phi)$, then by using Lemma \ref{him-p3-lem-1.13} we obtain the following Euclidean distance between $f(0)$ and the boundary of $f(\mathbb{D})$ as 
\begin{equation} \label{him-p3-e-3.13}
d(f(0), \partial f(\mathbb{D}))= \liminf \limits_{|z|\rightarrow 1} |f(z)-f(0)| \geq -h(-1).
\end{equation}
Since $f \in \mathcal{S}_{c}^{*}(\phi)$ and $\phi$ is starlike and symmetric with respect to real-axis, it follows that $g(z):=(f(z)+\overline{f(\bar{z})})/2$ belongs to $\mathcal{S}^{*}(\phi)$. Since $g \in \mathcal{S}^{*}(\phi)$, from Lemma \ref{him-p3-lem-1.6}, we have $g(z)/z \prec h(z)/z$. Therefore from Lemma \ref{him-p2-lem-1.23}, we obtain 
\begin{equation} \label{him-p3-e-3.8}
M_{g}(r) \leq M_{h}(r) \quad \mbox{for}\quad |z|=r\leq 1/3.
\end{equation}
From the definition of $\mathcal{S}_{c}^{*}(\phi)$, we have 
\begin{equation} \label{him-p3-e-3.9}
zf'(z)=g(z) \phi (\omega (z)),
\end{equation}
where $\omega:\mathbb{D} \rightarrow \mathbb{D}$ is analytic with $\omega(0)=0$. 
Since $\phi \circ \omega \prec \omega$, from Lemma \ref{him-p2-lem-1.23} we obtain
\begin{equation} \label{him-p3-e-3.10}
M_{\phi \circ \omega} (r) \leq M_{\phi} (r) \quad \mbox{for}\quad |z|=r\leq 1/3.
\end{equation}
A simplification of \eqref{him-p3-e-3.9} gives 
\begin{equation} \label{him-p3-e-3.11}
f(z)=\int \limits _{0}^{z} \dfrac{g(\xi)\phi(\omega(\xi))}{\xi}\, d\xi.
\end{equation}
Now, by making use of Lemma \ref{him-p3-lem-2.1} as well as \eqref{him-p3-e-3.8} and \eqref{him-p3-e-3.10} in \eqref{him-p3-e-3.11}, we obtain 
\begin{align} \label{him-p3-e-3.12}
|z|+\sum\limits_{n=2}^{\infty}|a_{n}||z|^{n}
& = M_{f}(r) \\ \nonumber
& \leq \int \limits _{0}^{r} \frac{M_{g}(t)M_{\phi \circ \omega}(t)}{t} \, dt \\ \nonumber
& \leq \int \limits _{0}^{r} \frac{M_{h}(t)M_{\phi}(t)}{t} \, dt 
\\ \nonumber
& =P(r)
\end{align}
for $|z|=r\leq 1/3$. We note that $P(r) \leq -h(-1)$, whenever $r \leq r_{f}$, where $r_{f}$ is the smallest positive root of $P(r)=-h(-1)$ in $(0,1)$. Going by the similar line of argument as in the proof of Theorem \ref{him-p3-thm-2.2}, the existance of the root $r_{f}$ is ensured by the following inequalities 
$$M_{h}(t)\geq |h(t)|, \, \, M_{h}(1) \geq |h(1)|\geq -h(-1) \quad \mbox{and} \quad M_{h}(0) < -h(-1).
$$ 
Thus, combining the inequalities \eqref{him-p3-e-3.12} and \eqref{him-p3-e-3.13} with the fact $P(r) \leq -h(-1)$ for $r \leq r_{f}$, we conclude that 
$$
|z|+\sum\limits_{n=2}^{\infty}|a_{n}||z|^{n} \leq d(f(0), \partial f(\mathbb{D}))
$$
for $|z|=r \leq \min \{1/3,r_{f}\}$.	
\end{pf}

\begin{pf} [{\bf Proof of Lemma  \ref{him-p3-lem-2.13}}]
From the definition of $\mathcal{S}_{cf}^{*}(\phi)$, we have $g \prec f$. Then by Lemma \ref{him-p2-lem-1.23}, we obtain $M_{g}(r) \leq M_{f}(r)$ for $|z|=r \leq 1/3$. Hence from \eqref{him-p3-e-2.11-a}, we obtain $\sum_{n=1}^{\infty}|g_{n}||z|^{n} \leq d(f(0), \partial f(\mathbb{D}))$ for $|z|=r \leq \min \{1/3, r_{f}\}$.
\end{pf}

\begin{pf} [{\bf Proof of Corollary \ref{him-p3-cor-2.15}}]
Since the coefficients of  $\phi (z)= (1+sz)^{2}$ with $0<s \leq 1/\sqrt{2}$ are all positive, in view of Remark \ref{him-p3-rem-2.2}, we obtain 
$$
P(r)=h(r)=r\exp \left(s\left(2r+\dfrac{sr^{2}}{2}\right)\right).
$$
Let $D_{1}(r)=h(r)+h(-1)$. Clearly $D_{1}$ is continuous in $[0,1]$. Observe that $D_{1}(0)<0$ and 
$$
D_{1}\left(\dfrac{1}{3}\right)=\dfrac{1}{3} \exp \left(s\left(\dfrac{s+12}{18}\right)\right)- \exp \left(s\left(-2+\dfrac{s}{2}\right)\right)>0,
$$
whenever $0.444981< s \leq 1/\sqrt{2}$. Therefore $D_{1}$ has a real root in $(0,1/3)$ and choose it to be $r_{f}$. Thus, from Remark \ref{him-p3-rem-2.2}, the radius $r_{f}$ is the best possible.
\end{pf}
\begin{pf} [{\bf Proof of Corollary \ref{him-p3-cor-2.17}}]
Let $\phi(z)=\alpha +(1-\alpha) e^{z}$. Then the coefficients of the Maclaurin series of $\phi(z)$ are positive for $0 \leq \alpha <1$. Let $D_{2}(r)=h(r)+h(-1)$, where 
$$
h\left(r\right)=r\exp\, \left((1-\alpha) \int\limits_{0}^{r}\left(\dfrac{-1+e^{t}}{t}\right)\, dt\right).
$$
It is easy to see that
$$
h\left(\dfrac{1}{3}\right)=\dfrac{1}{3} \exp\, \left((1-\alpha) \int\limits_{0}^{\frac{1}{3}}\left(\dfrac{-1+e^{t}}{t}\right)\, dt\right) \approx \dfrac{1}{3} (1.43807)^{1-\alpha}
$$
and 
$$
h(-1)=-\exp\, \left((1-\alpha) \int\limits_{0}^{-1}\left(\dfrac{-1+e^{t}}{t}\right)\, dt\right)\approx -(0.450859463)^{1-\alpha}.
$$
A simple computation shows that $D_{2}(1/3)=h(1/3)+h(-1)>0$ if $0 \leq \alpha < 0.05284$. Clearly, $D_{2}(0)=h(-1)<0$. 
Therefore, $D_{2}$ has a  root in $(0,1/3)$ and choose it to be $r_{f}$. In view of Remark \ref{him-p3-rem-2.2},  $r_{f}$ is the best possible.
\end{pf}

\begin{pf} [{\bf Proof of Corollary \ref{him-p3-cor-2.18}}]
Let $\phi(z)=\left((1+z)/(1-z)\right)^{\alpha}$ with $0< \alpha \leq 1$. From \cite{abu-2014}, it is guaranted that the coeffficients of the Maclaurin series of $\phi$ are positive. It is easy to see that 
$$
h(r)=r\exp \left(\int \limits _{0}^{r} \frac{\left(\frac{1+t}{1-t}\right)^{\alpha}-1}{t}\, dt \right).
$$
Then $D_{3}(r):=h(r)+h(-1)$ is continuous in $[0,1)$ and $D_{3}(0)<0$ and 
$D_{3}(1/3)=h(1/3)+h(-1)>0$. Thus $D_{3}$ has a  root in $(0,1)$ and choose it to be $r_{f}$. Hence, in view of Remark \ref{him-p3-rem-2.2},  $r_{f}$ is the best possible.
\end{pf}

\begin{pf} [{\bf Proof of Corollary \ref{him-p3-cor-2.19}}]
Let $\phi (z)=\left(1+(1-2\gamma)z\right)/(1-z)$. Then $h(z)=z/\left(1-z\right)^{2(1-\gamma)}$. It is easy to see that
$$h(1/3)= \dfrac{3^{2(1-\gamma)-1}}{2^{2(1-\gamma)}} \quad \mbox{and} \quad -h(-1)=\dfrac{1}{2^{2(1-\gamma)}}.
$$
Further, $h(1/3)>-h(-1)$ for $0 \leq \gamma \leq 1/2$. Therefore \eqref{him-p3-e-2.20} has a root in $(0,1/3)$ and monotonocity of $h$ ensures that this root is unique in $(0,1/3)$. Hence by the Remark \ref{him-p3-rem-2.2}, $r_{f}$ is the best possible for the class $\mathcal{S}_{c}^{*} \left(\left(1+(1-2\gamma)z\right)/(1-z)\right)$.
\end{pf}

\begin{pf} [{\bf Proof of Corollary \ref{him-p3-cor-2.21}}]
When $\phi (z)=(1+Az)/(1+Bz)$, then from \eqref{him-p2-e-3.8} we obtain 
\[ 
h(z)=
\begin{cases}
z(1+Bz)^{\frac{A-B}{B}}, \quad & B \neq 0\\[3mm]
ze^{Az}, &  B =0.
\end{cases}
\]
\begin{enumerate}
\item [(i)] When $B=0$, then $h(r)=re^{Ar}$ and $-h(-1)=e^{-A}$. We note that $h(1/3)>-h(-1)$ whenever $(1/3)e^{A/3}>e^{-A}$. That is when $A>(3/4)ln\, 3$. Therefore \eqref{him-p3-e-2.22} has a root in $(0,1/3)$ and choose $r_{f}$ be the smallest root in $(0,1/3)$. Hence $r_{f}$ is the best possible. 
\item [(ii)] If $B \neq 0$, then $h(r)=r(1+Br)^{(A-B)/B}$. It is easy to see that $h(1/3)>-h(-1)$ when $(1/3) \left(1+B/3\right)^{(A-B)/B} \geq \left(1-B\right)^{(A-B)/B}$. Therefore \eqref{him-p3-e-2.23} has a root in $(0,1/3)$ and choose $r_{f}$ to be the smallest root in $(0,1/3)$. Hence $r_{f}$ is the best possible.
\end{enumerate}	
\end{pf}

\begin{pf} [{\bf Proof of Theorem   \ref{him-p3-thm-2.24}}]
The proof of Theorem \ref{him-p3-thm-2.24} follows from  Theorem \ref{him-p3-thm-2.11} and the fact that $zf' \in \mathcal{S}_{c}^{*}(\phi)$ if, and only, if $f \in \mathcal{C}_{c}(\phi)$. For the bravity we complete the proof. Let $g(z):=(f(z)+\overline{f(\bar{z})})/2$. Since $\phi$ is starlike and symmetric with respect to real axis, $g \in \mathcal{C}(\phi)$. From the definition of $\mathcal{C}_{c}(\phi)$, we have 
\begin{equation} \label{him-p3-e-3.13-f}
\left(zf'(z)\right)'=g'(z)\phi (\omega (z)),
\end{equation}
where $\omega:\mathbb{D}\rightarrow \mathbb{D}$ is analytic with $\omega(0)=0$. A simple computation using \eqref{him-p3-e-3.13-f} shows that
\begin{equation} \label{him-p3-e-3.13-g}
f(z)=\int \limits _{0}^{z} \frac{1}{\xi} \int \limits _{0}^{\xi} g'(\eta) \phi(\omega(\eta))\, d\eta \, d\xi.
\end{equation}
Since $g\in \mathcal{C}(\phi)$, in view of Lemma \ref{him-p2-lem-1.13}, we have $g' \prec k'$ and hence by Lemma \ref{him-p2-lem-1.23}, we obtain 
\begin{equation} \label{him-p3-e-3.13-h}
M_{g'}(r)\leq M_{k'}(r) \quad \mbox{for } \quad r \leq 1/3.
\end{equation}
In view of Lemma \ref{him-p3-lem-2.1} and by using \eqref{him-p3-e-3.13-g} and \eqref{him-p3-e-3.13-h}, we obtain 
\begin{equation} \label{him-p3-e-3.13-i}
M_{f}(r) \leq \int \limits _{0}^{r}\frac{1}{s}\int \limits _{0}^{s}M_{k'}(t) M_{\phi}(t) \, dt\, ds=T(r) \quad \mbox{for } \quad r \leq 1/3.
\end{equation}
From Lemma \ref{him-p3-lem-1.14}, the Euclidean distance between $f(0)$ and the boundary of $f({\mathbb{D}})$ is
\begin{equation} \label{him-p3-e-3.13-j}
d(f(0), \partial f(\mathbb{D}))= \liminf \limits_{|z|\rightarrow 1} |f(z)-f(0)| \geq -k(-1).
\end{equation}
Clearly, $T(r) \leq -k(-1)$ for $r \leq r_{f}$, where $r_{f}$ is the smallest positive root of $T(r)=-k(-1)$ in $(0,1)$. Going by the similar lines of argument as in the proof of Theorem \ref{him-p3-thm-2.11}, the existance of the root $r_{f}$ is ensured by the following inequalities 
$$
M_{k}(r)\geq |k(r)|,\, \,  M_{k}(1) \geq |k(1)|\geq -k(-1) \quad \mbox{and} \quad M_{k}(0) < -k(-1).
$$
Therefore from \eqref{him-p3-e-3.13-i} and \eqref{him-p3-e-3.13-j}, we obtain 
$$
|z|+ \sum_{n=0}^{\infty} |a_{n}| |z|^{n}=M_{f}(r) \leq d(f(0),\partial f(\mathbb{D}))
$$
for $|z|=r\leq \min \{1/3,r_{f}\}$.
\end{pf}

\begin{pf} [{\bf Proof of Theorem   \ref{him-p3-thm-2.25}}]
Let $f \in \mathcal{C}_{s}(\phi)$, then it is evident that the Euclidean distance between $f(0)$ and the boundary of $f(\mathbb{D})$ is 
\begin{equation} \label{him-p3-e-3.13-e}
d(f(0), \partial f(\mathbb{D}))= \liminf \limits_{|z|\rightarrow 1} |f(z)-f(0)|\geq L_{s}(1).
\end{equation}
Since  $f \in \mathcal{C}_{s}(\phi)$ and $\phi$ is starlike and symmetric with respect to the real axis, then it follows that 
\begin{equation}
g(z):=\frac{f(z)-f(-z)}{2}=z+\sum \limits _{n=1}^{\infty} a_{2n+1}z^{2n+1} \, \in  \mathcal{C}(\phi).
\end{equation}
Here $g$ is an odd convex analytic function. Note that the function $K(z)=\int \limits _{0}^{z}(k'(\xi^{2}))^{1/2}\,\, d\xi$	defined in \eqref{him-p3-e-1.17} is an odd function in $\mathcal{C}(\phi)$. By Lemma \ref{him-p2-lem-1.13} we have $g' \prec K'$. Therefore from Lemma \ref{him-p2-lem-1.23}, we obtain
\begin{equation} \label{him-p3-e-3.13-a}
M_{g'}(r) \leq M_{K'}(r) \quad \mbox{for} \quad |z|=r\leq 1/3.
\end{equation}
From the definition of $\mathcal{C}_{s}(\phi)$, we have 
\begin{equation} \label{him-p3-e-3.13-b}
\left(zf'(z)\right)'=g'(z) \phi (\omega(z)).
\end{equation}
A simplication of \eqref{him-p3-e-3.13-b} gives 
\begin{equation} \label{him-p3-e-3.13-c}
f(z)=\int \limits _{0}^{z} \frac{1}{\xi} \int \limits_{0}^{\xi} g'(\eta) \phi(\omega(\eta))\, d\eta \, d\xi.
\end{equation}
By making use of Lemmas \ref{him-p3-lem-1.13} and \ref{him-p3-lem-2.1} and in view of \eqref{him-p3-e-3.13-a} and \eqref{him-p3-e-3.13-c}, we obtain 
\begin{align} \label{him-p3-e-3.13-d}
|z|+\sum\limits_{n=2}^{\infty}|a_{n}||z|^{n}=M_{f}(r)
 &\leq \int \limits _{0}^{r} \frac{1}{s} \int \limits_{0}^{s} M_{g'}(t) M_{\phi}(t) \,\, dt \,\, ds\\ \nonumber
 & \leq \int \limits _{0}^{r} \frac{1}{s} \int \limits_{0}^{s} M_{K'}(t) M_{\phi}(t) \,\, dt \,\, ds\\ \nonumber
 & = R_{s}(r) ,\nonumber
\end{align}
for $z|=r\leq 1/3$. A simple observation shows that $R_{s}(r) \leq L_{s}(1)$ for $r\leq r_{f}$, where $r_{f}$ is the smallest root of $R_{s}(r)=L_{s}(1)$ in $(0,1)$. The existance of the root is ensured by the following inequalities
 $$
M_{K'}(t)\geq |K'(t)|, \, \,R_{s}(1) \geq L_{s}(1) \quad \mbox{and} \quad R_{s}(0)\leq L_{s}(1)
$$ 
as well as the inequality \eqref{him-p3-e-13-a}. Using \eqref{him-p3-e-3.13-e} and \eqref{him-p3-e-3.13-d} with the fact that $R_{s}(r) \leq L_{s}(1)$ for $r\leq r_{f}$, we obtain 
$$
|z|+\sum\limits_{n=2}^{\infty}|a_{n}||z|^{n} \leq d(f(0), \partial f(\mathbb{D}))\quad \mbox{for} \quad |z|=r \leq \min \{1/3, r_{f}\}.
$$
This completes the proof.
\end{pf}

\noindent\textbf{Acknowledgement:}  The first author thanks SERB-MATRICS and the second author thanks CSIR for their support.

\end{document}